\documentclass[10pt]{amsproc}
\usepackage{amsfonts}
\usepackage{amsmath}
\usepackage{amssymb}
\usepackage{amsthm}
\usepackage{eucal}
\usepackage{graphicx}
\theoremstyle{plain}
\newtheorem{ansatz}{Ansatz Condition}
\newtheorem*{assumption}{Basic Assumption}
\newtheorem*{conjecture}{Conjecture}

\newtheorem{proposition}{Proposition}

\theoremstyle{definition}

\newcommand{\pd}{\partial}
\numberwithin{equation}{section}

\begin{document}
\title[Asymptotics for degenerate Ricci flow neckpinches]
{Formal matched asymptotics for\\degenerate Ricci flow neckpinches}
\author{Sigurd B.~Angenent}
\address[Sigurd Angenent]{University of Wisconsin-Madison}
\email{angenent@math.wisc.edu}
\urladdr{http://www.math.wisc.edu/\symbol{126}angenent/}
\author{James Isenberg}
\address[James Isenberg]{University of Oregon}
\email{isenberg@uoregon.edu}
\urladdr{http://www.uoregon.edu/\symbol{126}isenberg/}
\author{Dan Knopf}
\address[Dan Knopf]{University of Texas}
\email{danknopf@math.utexas.edu}
\urladdr{http://www.ma.utexas.edu/users/danknopf}
\thanks{S.B.A.~acknowledges NSF support via DMS-0705431.
        J.I.~ acknowledges NSF support via PHY-0652903 and PHY-0968612.
        D.K.~ acknowledges NSF support via DMS-0545984.}
        %%\emph{This is work in progress and is not for distribution.} (@~\the\time~on \today)

\begin{abstract}
Gu and Zhu \cite{GZ08} have shown that Type-II Ricci flow singularities develop from
nongeneric rotationally symmetric Riemannian metrics on $\mathcal{S}^{n+1}\,(n\geq 2)$.
In this paper, we describe and provide plausibility arguments for a detailed asymptotic
profile and rate of curvature blow-up that we predict such solutions exhibit.
\end{abstract}

\maketitle
\tableofcontents

\section{Introduction}
\label{Introduction}

Let $(\mathcal{S}^{n+1},g(t):0\leq t<T)$ be a rotationally symmetric
solution of Ricci flow. Gu and Zhu prove that Type-II singularities can
occur for special nongeneric initial data of this type \cite{GZ08}.
Garfinkle and one of the authors provide numerical simulations of the
formation of such singularities at either pole \cite{GI03, GI07}. However,
almost nothing is known about the asymptotics of such singularity formation,
except for complete noncompact solutions on $\mathbb{R}^2$, where Ricci flow
coincides with logarithmic fast diffusion, $u_t=\Delta \log u$. (The asymptotics
of logarithmic fast diffusion, which are unrelated to the results in this paper,
were derived by King \cite{King93} and subsequently proved in $\mathbb{R}^2$ by
Daskalopoulos and \v{S}e\v{s}um \cite{DS10}.)

We say a sequence $\{(p_i,t_i)\}_{i=0}^\infty$ of points and times in a Ricci flow
solution is a \emph{blow-up sequence} at time $T$ if $t_i \nearrow T$ and
$|\mathrm{Rm}(x_i,t_i)|\rightarrow\infty$ as $i\rightarrow\infty$. We say
$(\mathcal{M}^{n+1},g(t))$ develops a \emph{neckpinch singularity} at $T<\infty$ if there
is some blow-up sequence at $T$ whose corresponding sequence of parabolic dilations has as
a pointed limit the self-similar Ricci soliton on the cylinder $\mathbb{R}\times\mathcal{S}^n$.
We call a neckpinch \emph{nondegenerate} if every complete pointed limit formed from a blowup
sequence at $T$ is a solution on the cylinder; we call it \emph{degenerate} if a complete smooth
limit of some blow-up sequence has another topology. Rotationally symmetric nondegenerate
neckpinches have been studied by Simon \cite{Simon00} and two of the authors \cite{AK04, AK07}, who
proved they are Type-I (``rapidly forming'') singularities in which the curvature blows up at
the natural parabolic rate, $(T-t)\sup_{p\in\mathcal{M}^{n+1}}|\mathrm{Rm}(p,t)|<\infty$.
On the other hand, Type-II (``slowly forming'') singularities have the property that
$(T-t)\sup_{p\in\mathcal{M}^{n+1}}|\mathrm{Rm}(p,t)|=\infty$. The compact Type-II singularities
proved to exist by Gu and Zhu \cite{GZ08} are degenerate neckpinches.

This paper is the first of two in which we study the formation of degenerate Ricci flow
neckpinch singularities. In the present work, we assume that a degenerate neckpinch
singularity occurs at a pole at time $T<\infty$, and we derive formal matched asymptotics
for the solution as it approaches the singularity. This procedure provides evidence for
a conjectural picture of the behavior of some (not necessarily all) solutions that develop
rotationally symmetric degenerate neckpinch singularities. In particular, it predicts
precise rates of Type-II curvature blow-up.\footnote
{For a comprehensive statement of our results, see Section~\ref{TheBitterFinalEnd}.}
In forthcoming work, we will provide rigorous proof that there exist solutions exhibiting
the asymptotic behavior formally described here.

Given any rotationally symmetric family $g(t)$ of metrics on $\mathcal{S}^{n+1}$,
one may remove the poles $P_{\pm}$ and write the metrics on
$\mathcal{S}^{n+1}\backslash\{P_{\pm}\}\approx(-1,1)\times\mathcal{S}^{n}$ in the form
$g(t)=(\mathrm{d}s)^{2}+\psi^{2}(s,t)\,g_{\rm{can}}$, where $s(x,t)$ denotes
$g(t)$-arclength to $x\in [-1,1]$ from a fixed point $x_0 \in (-1,1)$. (Here, $g_{\rm{can}}$
is the canonical unit sphere metric on $\mathcal{S}^n$; see Section~\ref{Notation}
for a detailed discussion of these coordinates.) Thus the function $\psi(s(x,t),t)$
completely characterizes a given solution. Our Basic~Assumption in this paper,
explained in detail in Section~\ref{Notation}, is that the initial metric $g(s,0)$,
hence $\psi(s,0)$, satisfies certain curvature restrictions which ensure that the geometries
we consider are sufficiently close to those studied in \cite{GI03, GI07} and \cite{GZ08},
and therefore are likely to develop neckpinches. It also allows us to employ certain prior
results of two of the authors \cite{AK04}, which are useful for the arguments made
here. (Compare \cite{AK07}.)

In the language of Section~\ref{Notation}, let $\hat{s}(t)$ denote the location of the
local maximum (``bump'') of $\psi$ closest to the right pole $x=1$. Note that
$\hat{s}(t)$ may be only an upper-semicontinuous function of time, because a bump
and an adjacent neck can join and annihilate each other. Lemma~7.1 of \cite{AK04}
proves that the solution exists until $\psi$ becomes zero somewhere other than at the
poles. Since there is always at least one positive local maximum of $\psi$, the quantity
$\hat{s}(t)$ is defined for as long as the solution exists. Lemma~7.2 of \cite{AK04}
proves that if $\lim_{t\nearrow T}\psi(\hat{s}(t),t)>0$, then no singularity occurs at the
right pole. Hence, if a degenerate neckpinch does happen at the right pole, it must be that
$\lim_{t\nearrow T}\psi(\hat{s}(t),t)=0$. This can happen only if (i) the radius $\psi$ vanishes
on an open set, or (ii) the bump marked by $\hat{s}$ moves to the right pole. Note that these
alternatives are not mutually exclusive. In either case, it follows that there are points
$\tilde{s}(t)$ such that
%% \footnote{Here and throughout this paper, partial derivatives are denoted by subscripts:
%% i.e.,  $\psi _{s}(s,t)$ denotes $\frac{\partial \psi}{\partial s}(s,t)$.}
\begin{equation} \label{javelin}
\lim_{t\nearrow T}\tilde{s}(t)=s(1,T)
\qquad \text{and} \qquad
\lim_{t\nearrow T}\psi_{s}(\tilde{s}(t),t)=0.
\end{equation}
Observe that \eqref{javelin} is incompatible with the boundary condition
\eqref{psi-s-boundary} necessary for regularity at the pole, which forces
$\psi=(s(1)-s)+o\big(s(1)-s\big)$ as $s\nearrow s(1)$. Therefore, a consequence
of our Basic~Assumption is that if a degenerate singularity does develop, then this
expansion for $\psi$ cannot hold uniformly in $s$ as $t\nearrow T$. Instead, the solution
must exhibit different qualitative behaviors in a sequence of time-dependent spatial regions.

Starting with that fact, we construct in this paper a conjectural model for
rotationally symmetric Ricci flow solutions that develop a degenerate neckpinch,
and we check its consistency. We do this systematically by studying approximate asymptotic
expansions to the solution in four connected regions, which we call (moving in from the
pole) the \emph{tip}, \emph{parabolic}, \emph{intermediate,} and \emph{outer} regions.
By matching these expansions at the intersections of the sequential regions, we produce our
model. This is the essence of the formal matched asymptotics process familiar to
applied mathematicians. The process involves first formulating an \emph{Ansatz}, which consists of
a series of assumptions (listed below) pertaining to the geometry and its evolution equations in the
four regions. Then, matching across the boundaries of the regions, one constructs formal (approximate)
solutions. Finally, one checks that these formal solutions remain consistent with the \emph{Ansatz},
and argues that the approximations remain sufficiently accurate.
\medskip

For clarity of exposition, we establish notation in Section~\ref{Notation} and then begin our study
working in the parabolic region rather than the tip. We treat the tip (the most critical region) last.
Finally, in Section~\ref{TheBitterFinalEnd}, we provide a detailed summary of our results and the
conjectural picture they provide.

\section{Basic equations}
\label{Notation}

We begin by recalling some basic identities for $\rm{SO}(n+1)$-invariant Ricci flow solutions.

To avoid working in multiple patches, it is convenient to puncture the sphere
$\mathcal{S}^{n+1}$ at the poles $P_{\pm}$ and identify
$\mathcal{S}^{n+1}\backslash\{P_{\pm}\}$ with $(-1,1)\times\mathcal{S}^{n}$.
If we let $x$ denote the coordinate on the interval $(-1,1)$ and let
$g_{\rm{can}}$ denote the canonical unit sphere metric, then an arbitrary family $g(t)$ of
smooth $\rm{SO}(n+1)$-invariant metrics on $\mathcal{S}^{n+1}$ may be written in geodesic
polar coordinates as
\begin{equation} \label{unnatural-metric}
g(t)=\varphi^{2}(x,t)\,(\mathrm{d}x)^{2}+\psi^{2}(x,t)\,g_{\rm{can}}.
\end{equation}
Denoting the distance from the equator $\{0\}\times\mathcal{S}^{n}$ by
\[
s(x,t)=\int_{0}^{x}\varphi(\xi,t)\,\rm{d}\xi ,
\]
allows one to write \eqref{unnatural-metric} in the more geometrically natural form
\begin{equation}
g=(\mathrm{d}s)^2 +\psi^{2}(s(x,t),t)\,g_{\rm{can}},  \label{natural-metric}
\end{equation}
where ``$\mathrm{d}$'' is the differential with respect to the space variables
but not the time variables.  We also agree to write
\[
\frac{\pd}{\pd s} = \frac{1}{\varphi(x,t)}\frac{\pd}{\pd x}.
\]
With this convention, $\pd/\pd s$ and $\pd/\pd t$ do not commute; instead one gets
the usual commutator
\[
\left[ \frac{\pd}{\pd t}, \frac{\pd}{\pd s} \right] =
-\frac{\varphi_t}{\varphi} \frac{\pd}{\pd s}.
\]
Smoothness at the poles requires that $\psi$ satisfy the boundary conditions
\begin{equation}
\left. \psi _{s}\right\vert_{x=\pm 1}=\mp 1.  \label{psi-s-boundary}
\end{equation}

The metric \eqref{natural-metric} has two distinguished sectional curvatures:
the curvature
\begin{equation}
L=\frac{1-\psi_{s}^{2}}{\psi^{2}}
\end{equation}
of a plane tangent to $\{s\}\times \mathcal{S}^{n}$, and the curvature
\begin{equation}
K=-\frac{\psi _{ss}}{\psi}
\end{equation}
of an orthogonal plane.

Under Ricci flow,
\[
\frac{\partial}{\partial t}g=-2\rm{Rc},
\]
the quantities $\varphi$ and $\psi$ evolve by the degenerate parabolic system
\begin{align*}
\varphi_t &=    n\left(\frac{\psi_{xx}}{\varphi\psi}
                -\frac{\varphi_x \psi_x}{\varphi^2 \psi}\right),\\
\psi_t &=   \frac{\psi_{xx}}{\varphi^2}
            -\frac{\varphi_x \psi_x}{\varphi^3}+(n-1)\frac{\psi_x^2}{\varphi^2 \psi}-\frac{n-1}{\psi}.
\end{align*}
The degeneracy is due to invariance of the system under the infinite-dimensional
diffeomorphism group. By writing the spatial derivatives in terms of $\pd/\pd s$,
one can simplify the appearance of these equations, effectively by fixing arclength
as a gauge. The quantities $\varphi$ and $\psi$
in \eqref{unnatural-metric} then evolve by
\begin{subequations}
  \begin{align}
    \varphi_{t}&=n\frac{\psi_{ss}}{\psi}\varphi,
    \label{phi-s-evolution}\\
    \psi_{t}&=\psi_{ss}-(n-1)\frac{1-\psi_{s}^{2}}{\psi},
    \label{psi-s-evolution}
  \end{align}
\end{subequations}
respectively.
\medskip

In \cite{AK04}, nondegenerate neckpinch singularity formation is established for an open
set of initial data of the form \eqref{natural-metric} on $\mathcal{S}^{n+1}$
$(n\geq 2)$ satisfying the following assumptions:\footnote
    {We call local minima of $\psi$ ``necks'' and local maxima ``bumps''.
    The ``polar caps'' are the regions on either side of the outermost bumps.}
\begin{enumerate}
\item The sectional curvature $L$ of planes tangent to each sphere
  $\{s\}\times\mathcal{S}^{n}$ is positive.

\item The Ricci curvature
    $\mathrm{Rc}=nK(\mathrm{d}s)^2 +[K+(n-1)L]\psi^2\,g_{\mathrm{can}}$
    is positive on each polar cap.

\item The scalar curvature $R=2nK+n(n-1)L$ is positive everywhere.

\item The metric has at least one neck and is ``sufficiently pinched''
  in the sense that the value of $\psi$ at the smallest neck is
  sufficiently small relative to its value at either adjacent bump.

\hspace*{-1.5 cm} In \cite{AK07}, precise asymptotics are derived under the additional
  assumption:

\item The metric is reflection symmetric, $\psi(s)=\psi(-s)$,
  and the smallest neck is at $x=0$.
\end{enumerate}
\medskip

For our analysis here of degenerate neckpinch solutions, we start from initial geometric data
that satisfy some but not all of these assumptions.  In particular, we drop the condition on the
``tightness'' of the initial neck  pinching, and we also do not require reflection symmetry.\footnote
    {As seen in Section~\ref{HermiteAgain}, some though not all of the solutions we construct
    here are reflection symmetric.}
We now state our underlying assumptions on the initial geometry:

\begin{assumption} \label{BasicAssumptions}
The solutions of Ricci flow we consider in this paper are $\mathrm{SO}(n+1)$-invariant.
They satisfy conditions~(1)--(3) above, and they have initial data with at least one neck.
Furthermore, we assume that a singularity occurs at the right pole $(x=+1)$ at some time $T<\infty$.
\end{assumption}

\section{The parabolic region} \label{HermiteAgain}

A consequence of our Basic~Assumption is that $R>0$ at $t=0$. By standard
arguments, this implies that the solution becomes singular at some time $T<\infty$.

We assume that a neckpinch occurs at some $x_0$; and by diffeomorphism invariance, we
may assume that $x_0=0$. Thus in the parabolic region to be characterized below,
it is natural to study the system in the renormalized variables $\tau$ and $\sigma$
defined by
\begin{equation}
\tau(t) :=-\log(T-t),
\end{equation}
and
\begin{equation}
\label{define-sigma}
\sigma(x,t) :=\frac{s(x,t)}{\sqrt{T-t}}=e^{\tau/2}s(x,t),
\end{equation}
where $s(x,t)$ denotes arclength from the point $x=0$.
Note that $\sigma$ is a natural choice for the parabolic region, in
which time scales like distance squared. Use of $\tau$ is not necessary
but is convenient for the calculations that follow. (Compare \cite{AK07}.)

The evolution equation satisfied by $\psi$, equation~\eqref{psi-s-evolution},
implies that at the local maximum (``bump'') $\hat s(t)$ closest to the right pole,
one has $\psi_{t}(\hat{s},t)\leq(n-1)/\psi $. On the other hand, applying a version
of the maximum principle to the same equation shows that the radius of the smallest
neck satisfies the condition $\psi_{\min}(t)\geq\sqrt{(n-1)(T-t)}$.
(For the proof, see Lemma~6.1 of \cite{AK04}.)
Hence, we introduce the rescaled radius
\begin{equation}
\label{define-u}
U(\sigma,\tau):=\frac{\psi(s,t)}{\sqrt{2(n-1)(T-t)}}.
\end{equation}
Note that $U$ is a function of the coordinate $\sigma$. In Appendix~\ref{WhyTheHeck?},
we compute the evolution equation satisfied by $U$, which is
\begin{equation}
  \label{U-evolution}
  U_{\tau}
  =
  U_{\sigma\sigma} -\left(\frac{\sigma}{2} + nI\right) U_{\sigma}
  + (n-1)\frac{U_{\sigma}^{2}}{U} + \tfrac12 \left(U- \frac{1}{U}\right),
\end{equation}
where
\[
I := \int_0^\sigma \frac{U_{\sigma\sigma}}{U}\,\mathrm{d}\sigma.
\]
Note that equation~\eqref{U-evolution} is parabolic, but contains a nonlocal term.
\medskip

Motivated by Perelman's work \cite{P1} in dimension $n+1=3$, one expects the solution
to be approximately cylindrical a controlled distance away from the pole, depending
only on the curvature there. (Essentially, this is contained in Corollary~11.8 of \cite{P1};
for more details and a refinement, see Sections~47--48 and in particular Corollary~48.1 of
\cite{KL08}.) This expectation effectively determines the behavior we anticipate in
the parabolic region $\Omega_{\rm{par}}$.

\begin{ansatz}
\label{MainParabolic}
As $\tau\nearrow\infty$ (equivalently, as $t\nearrow T$), the rescaled radius
$U$ converges uniformly in regions $|\sigma|\leq \textrm{const}$, where
  \[
  U(\sigma, \tau) \to 1.
  \]
\end{ansatz}

Since an exact cylinder solution is given by $U(\sigma,\tau)=1$, the parabolic region
$\Omega_{\rm{par}}$ may be characterized as that within which the quantity
\begin{equation}
V(\sigma,\tau):=U(\sigma,\tau)-1  \label{define-v}
\end{equation}
is small.

By substituting $U=1+V$ in equation~\eqref{U-evolution},\footnote{See equation~\eqref{v-evolution}.}
one finds that $V$ evolves by
\begin{equation}
V_{\tau}=AV+N(V),  \label{v-tau-parabolic}
\end{equation}
where $A$ is the linear operator
\begin{equation}
AV := V_{\sigma \sigma }-\frac{\sigma }{2}V_{\sigma }+V,
\label{Linear-Parabolic}
\end{equation}
and $N$ consists of nonlinear terms,
\begin{equation}
N(V):=\frac{2(n-1)V_{\sigma}^2-V^2}{2(1+V)}-nIV_\sigma,  \label{Nonlinear-Parabolic}
\end{equation}
with $I = \int_{0}^\sigma V_{\sigma\sigma}/(1+V)\,\rm{d}\sigma$ its sole
non-local term. %%\footnote{See equation~\eqref{non-local}.}

The operator $A$ appears in many blow-up problems for semilinear
parabolic equations of the type $u_t=u_{xx}+ u^p$ (e.g.~\cite{GK85}),
and also in the analysis of mean curvature flow neckpinches (e.g.~\cite{AV97}).
The operator $A$ is self adjoint in $L^{2}(\mathbb{R},\,e^{-\sigma ^{2}/4}\,\mathrm{d}\sigma)$.
It has pure point spectrum with eigenvalues $\{\lambda_k\}_{k=0}^{\infty}$, where
\begin{equation}
\lambda_{k}:=1-\frac{k}{2}.
\end{equation}
The associated eigenfunctions are the Hermite polynomials $h_{k}(\sigma)$.
We adopt the normalization that the highest-order term has coefficient $1$,
so that $h_{0}(\sigma)=1$, $h_{1}(\sigma)=\sigma$, $h_{2}(\sigma)=\sigma^2-2$,
$h_3(\sigma)=\sigma^3 - 6\sigma$, and in general,
\[
h_{k}(\sigma )=\sum_{j=0}^{k}\eta_{j}\sigma ^{j},
\]
for certain determined constants $\eta_j$, with $\eta_{k}=1$. Recall that this
sum for the Hermite polynomial $h_k$ contains only odd powers of $\sigma$ if
$k $ is odd, and even powers if $k$ is even.

In the parabolic region $\Omega_{\rm{par}}$, which may be regarded as the space-time
region where one expects $V\approx 0$, the linear term $AV$ should dominate
the nonlinear term $N(V)$, as long as $V$ is orthogonal to the kernel of $A$.
(We argue below that this orthogonality follows from the \emph{Ansatz} adopted
in this paper.) Therefore, we may assume that $V\approx \tilde{V}$, where
$\tilde{V}$ is an exact solution of the approximate (linearized) equation
$\tilde{V}_{\tau}=A\tilde{V}$.

The general solution to the linear \textsc{pde} $ \tilde{V}_{\tau }=A\tilde{V}$ can be written
in the form of a Fourier series by expanding in the Hermite eigenfunctions, namely
$\tilde{V}(\sigma,\tau)=\sum_{k=0}^{\infty}b_k e^{\lambda_k \tau}h_k (\sigma)$,
where all coefficients $b_k$ except possibly $b_2$ are constants. If we treat $\tilde{V}$ as the
first term in an approximate expansion for solutions of the nonlinear system~\eqref{v-tau-parabolic},
then we are led to write
\begin{equation}
\label{expansion}
V(\sigma,\tau)\approx b_0 e^{\tau}+b_1 e^{\tau /2}h_1 (\sigma)+b_2 (\tau)h_2 (\sigma)
    +\sum_{k=3}^{\infty}b_k e^{\lambda_k \tau}h_k (\sigma),
\end{equation}
where $b_2 (\tau)$ is allowed to be a function of time, since it corresponds to
a null eigenvalue of the linearized operator and hence to motion on a center
manifold where the contributions of $N(V)$ cannot be ignored.\footnote{Recall
that $b_2 (\tau)=(8\tau)^{-1}$ for the rotationally symmetric neckpinch \cite{AK07}.}
In a standard initial value problem, the coefficients $b_k$ are determined by the initial
data. Here, they may be determined by matching at the intersection with the adjacent
regions.

We now make a pair of additional assumptions.
\begin{ansatz}
\label{Parabolic-Ansatz}
(i) The solution does not approach a cylinder too quickly, in the sense that
 $\sup_{(\sigma,\tau)\in\Omega_{\mathrm{par}}}|V(\sigma,\tau)|\geq e^{-C\tau}$ for some $C<\infty$.
(ii) The singularity time $T<\infty$ depends continuously on the initial data $g_0$.
\end{ansatz}

The key consequence of Ansatz Condition~\ref{Parabolic-Ansatz}, which we explain below,
 is that one eigenmode dominates in the sense that
\begin{equation} \label{eigenmode}
V(\sigma ,\tau )\approx V_k (\sigma,\tau):=b_k e^{\lambda_k \tau}h_k (\sigma).
\end{equation}
Because $A$ has pure point spectrum in $L^{2}(\mathbb{R},\,e^{-\sigma^{2}/4}\,\mathrm{d}\sigma)$,
this is reasonable and consistent with Ansatz Condition~\ref{MainParabolic} and with the fact,
proved in \cite{AK04}, that singularity formation at the pole is nongeneric in the class of all
rotationally symmetric solutions.

In fact, we may assume that equation~\eqref{eigenmode} holds for some $k\geq 3$.
Here is why. Part (i) of Ansatz Condition~\ref{Parabolic-Ansatz} implies that
some $b_k\neq 0$. (Compare Section~2.16 of \cite{AK07}.) Part (ii) implies that one
can without loss of generality choose the free parameter $T=T(g_0)$ so that $b_0=0$.
(This too is true for the nondegenerate neckpinch; see Section~2.15 of \cite{AK07}.)
Ansatz Condition~\ref{MainParabolic} implies that $b_1=0$. In Section~\ref{Tip},
we discard $k=2$ as a consequence of Ansatz Condition~\ref{TipAnsatz}. Hence we have
$k\geq 3$, which implies orthogonality to the kernel of $A$, as promised above.

We can now obtain a rough estimate for the size of the parabolic region. Recall that
the main assumption made for the parabolic region, Ansatz Condition~\ref{MainParabolic},
applies only so long as $U \approx 1$, hence as long as
\begin{equation} \label{Parabolic-Intermediate-Interface}
V(\sigma,\tau)= b_k e^{\lambda_k \tau}\sigma^k +
\mathcal{O}(e^{\lambda_k \tau}\sigma^{k-2})\qquad(\sigma\rightarrow\infty)
\end{equation}
is small, hence only as long as $|\sigma|\ll e^{(\frac{1}{2}-\frac{1}{k})\tau}$.
We label the region characterized by larger values of $|\sigma|$ the intermediate region.

\section{The intermediate region}

The assumptions we have made to obtain a formal expansion in the parabolic
region are expected to be valid only so long as $V(\sigma,\tau)$ is much smaller
than one. We thus leave the parabolic region and enter the next region
for those values of $(\sigma, \tau)$  for which $V(\sigma,\tau)$ is of order one.

By Ansatz Condition~\ref{Parabolic-Ansatz} and its consequence \eqref{eigenmode},
the dominant term in the Fourier expansion of $V$ is a polynomial with leading term
comparable to $e^{\lambda _{k}\tau }\sigma ^{k}$. So we find it useful to introduce a new
space variable
\begin{equation} \label{define-rho}
\rho :=(e^{\lambda_k\tau})^{1/k}\sigma=e^{(\frac{1}{k}-\frac{1}{2})\tau}\sigma,
\end{equation}
and use $\rho$ to demarcate the intermediate ($\rho \approx 1$) region.\

In the intermediate region $\Omega_{\rm{int}}$, we define
\begin{equation} \label{define-W}
W(\rho,t)=U(\sigma,\tau)\quad \big(=1+V(\sigma,\tau)\big).
\end{equation}
Recalling the evolution equation~\eqref{U-evolution} satisfied by $U$, one computes
that $W$ satisfies the \textsc{pde}
\begin{multline*}
\frac{1}{2}(W-W^{-1})-\frac{\rho}{k} W_{\rho}\\
    =e^{-\tau}W_t
    -e^{(2/k-1)\tau}
    \left\{W_{\rho\rho}+(n-1)\frac{W_{\rho}^{2}}{W}-n\int_0^\rho\frac{W_{\rho\rho}}{W}\,\rm{d}\rho\right\}.
\end{multline*}
This appears to be a small time-dependent perturbation of an \textsc{ode} in
$\rho$. We make this intuition precise in the following additional assumption:

\begin{ansatz}
\label{IntermediateAnsatz}
$W_t$, $W_{\rho}$, and $W_{\rho\rho}$ are all bounded in the intermediate region.
\end{ansatz}

Ansatz Condition~\ref{IntermediateAnsatz} implies that
$\frac{\rho}{k}W_{\rho}-\frac{1}{2}W+\frac{1}{2}W^{-1}=o(1)$ as $\tau\rightarrow\infty$.
This suggests that $W(\rho,t)\approx \tilde{W}(\rho,t)$,
where $\tilde{W}$ solves the first-order \textsc{ode}
\begin{equation} \label{w-approximate}
\frac{\rho}{k}\tilde{W}_{\rho}-\frac{1}{2}\tilde{W}+\frac{1}{2}\tilde{W}^{-1}=0.
\end{equation}
One readily determines that the general solution to equation~\eqref{w-approximate} is
\begin{equation}
\tilde{W}(\rho, t )=\sqrt{1-(\rho/c)^{k}}.
\end{equation}
Notice that in solving this \textsc{ode}, we obtain a  ``constant of integration''
$c(t)$ which may \emph{a priori} depend on time. Matching considerations will show
that $c$ is in fact constant. (The minus sign above is chosen so that the singularity
occurs at the right pole if $c>0$.)

In order to match the fields in the intermediate region with those in the
parabolic region, it is useful to consider the asymptotic expansion of
$\tilde{W}$ about $\rho=0$; one gets
\[
\tilde{W}(\rho,t)=1-\frac{1}{2}(\rho/c(t))^{k}-\frac{1}{8}(\rho/c(t))^{2k}-\cdots.
\]
Thus for $0<\rho=(e^{\lambda_{k}\tau})^{1/k}\sigma\ll 1$,
Ansatz Condition~\ref{IntermediateAnsatz} implies that
\begin{equation}
\label{W-approximation}
W(\rho,t)\approx
    \tilde{W}(\rho,t)=1-\frac{c(t)^k}{2}e^{\lambda_k\tau}\sigma^k +\cdots,
\end{equation}
which matches the parabolic expansion~\eqref{Parabolic-Intermediate-Interface}
if $c(t)$ is determined by $b_k$, which is constant for all $k\geq 3$. It is convenient
in what follows to regard $b_k$ as determined by $c>0$; thus one has
\begin{equation} \label{choose-b}
b_k = -\frac{1}{2}c^{-k}.
\end{equation}

Going in the other direction, we expect to transition to the tip region where $W\approx 0$,
namely where $\rho\approx c$. It is unsurprising that one gets the same conclusion from
equation~\eqref{Parabolic-Intermediate-Interface} by solving $V\approx-1$ for large $|\sigma|$.

\section{The outer region}

For $\rho\gg 1$, which marks the outer region $\Omega_{\rm{out}}$, we unwrap definitions
\eqref{define-sigma}, \eqref{define-u}, \eqref{define-rho}, and \eqref{define-W} to obtain
\begin{equation} \label{DeathProfile}
\psi(s,t)^2
    =2(n-1)e^{-\tau}W(e^{\tau/k}s,\tau)^2
    \approx2(n-1)\left[(T-t)-\left(\frac{s}{c}\right)^k\right].
\end{equation}

If $k\geq3$ is odd, then equation~\eqref{DeathProfile} implies that the $g(t)$-measure
of the open set $(0,1)\times\mathcal{S}^{n}$ goes to zero as $t\nearrow T$.
Moreover, the $t=T$ limit profile of $\psi$ can vanish to arbitrarily high order, behaving
like $s^{k/2}$ as $s\nearrow 0$.

If $k\geq4$ is even, then equation~\eqref{DeathProfile} implies that the solution encounters
a global singularity at $t=T$, with the whole manifold shrinking into a non-round point.

In either case, a set of nonzero $g(0)$-measure is destroyed, in contrast to the nondegenerate
Type-I neckpinches studied in \cite{AK04, AK07}, which become singular only on the hypersurface
$\{x_0\}\times\mathcal{S}^n$. Note that it is possible for local Type-I singularities to
destroy sets of nonzero measure. See examples by one of the authors \cite{FIK03}, as well as
recent work of Enders, M\"{u}ller, and Topping \cite{EMT10}.

\section{The tip region}\label{Tip}

Obtaining precise asymptotics at the right pole is complicated by the fact that the natural
geodesic polar coordinate system~\eqref{unnatural-metric} becomes singular there. As in
\cite{ACK09}, we overcome this difficulty by choosing new local coordinates.

The fact that $\psi_{s}(s(1,t),t)=-1$ for all $t<T$ implies that $\psi _{s}<0$ in a
small time-dependent neighborhood of the right pole $x=1$. So we may regard $\psi(s,t)$
as a new local radial coordinate, thereby regarding $s$ as a function of $\psi$ and $t$.
(Compare \cite{AV97} and \cite{ACK09}.) More precisely, there is a function
$y(\psi,t)<0$ defined for small $\psi>0$ and times $t$ near $T$ such that
\begin{equation}
\psi_{s}(s,t)=y(\psi(s,t),t).
\end{equation}
In terms of this coordinate, the metric takes the form
\begin{equation}
g=y(\psi(s,t),t)^{-2}(\mathrm{d}\psi)^{2}+\psi^{2}\,g_{\rm{can}},
\end{equation}
with $y(\psi,t)$ the unknown function whose evolution controls the geometry near the pole.

For convenience, we replace $y<0$ by the quantity $z=y^{2}$, which evolves by the \textsc{pde}
$z_{t}=\mathcal{F}_{\psi }[z]$, where
\begin{equation}   \label{Define-F}
\mathcal{F}_{\psi }[z]
    :=\frac{1}{\psi^2}\left\{\psi^2 z z_{\psi\psi }
    -\frac{1}{2}(\psi z_\psi)^2 + (n-1-z)\psi z_{\psi}
    +2(n-1)(1-z)z\right\}.
\end{equation}

For the tip region, we now introduce a $t$-dependent expansion factor for the radial coordinate, setting
\begin{equation}
\gamma(s,t):=\Gamma(\tau(t))\,\psi(s,t),
\end{equation}
with  $\Gamma$ to be determined below by matching considerations.
Defining
\begin{equation}
Z(\gamma,t):=z(\psi,t)\quad \big(=y^2 (\psi,t)\big),
\end{equation}
one determines from equation~\eqref{Define-F} that $Z$ satisfies
\begin{equation} \label{Z-time-evolution}
\Gamma^{-2}(Z_{t}+\Gamma^{-1}\Gamma_\tau\gamma Z_{\gamma})
=\mathcal{F}_{\gamma}[Z],
\end{equation}
where $\mathcal{F}_{\gamma}[\cdot]$ is the operator appearing in equation~\eqref{Define-F},
with all $\psi$ derivatives replaced by $\gamma$ derivatives.

As we observe in Appendix~\ref{Bryant}, the \textsc{ode} $\mathcal{F}_{\gamma}[\tilde{Z}]=0$
admits a one-parameter family of complete solutions satisfying the boundary conditions
$\tilde{Z}(0)=1$ and $\tilde{Z}(\infty)=0$. Any such solution is given by
\[
\tilde{Z}(\gamma)=B\left(\frac{\gamma}{a}\right),
\]
where $a>0$ is an arbitrary
scaling parameter, and $B$ is the profile function of the \emph{Bryant steady soliton} metric,
\[
g_{\mathrm{B}} = B(r)^{-1}(\mathrm{d}r)^{2}+r^{2}\,g_{\rm{can}},
\]
discovered in unpublished work of Bryant, who proved that it is, up to homothety, the
unique complete rotationally symmetric non-flat steady gradient soliton on $\mathbb{R}^{n+1}$
for $n \geq 2$. Recent results of Cao and Chen (followed by a simplified proof by Bryant) allow
one to replace the words ``rotationally symmetric'' in the uniqueness statement above by the
words ``locally conformally flat'' \cite{CC10}. Uniqueness under the assumption of local
conformal flatness for $n+1\geq4$ was also proved independently by Catino and Mantegazza.
Recent progress by Brendle  allows one to replace local conformal flatness by the condition
that a certain vector field $V:=\nabla R+\varrho(R)\nabla f$ decays fast enough at spatial
infinity \cite{Brendle11}. (Here, $f$ is the soliton potential function, and $\varrho(R)$
is chosen so that $V$ vanishes on the Bryant soliton.) We briefly review some relevant
properties of the Bryant soliton in Appendix~\ref{Bryant}.

Numerical studies for Ricci flow of rotationally symmetric neckpinches \cite{GI03, GI07}
strongly support the contention that for certain initial data, the flow near the tip approaches
a Bryant soliton model. Results of Gu and Zhu also support this expectation \cite{GZ08}.
Therefore, in the tip region, we adopt the following assumptions, which suggest that a solution
of \eqref{Z-time-evolution} should be approximated by a steady-state solution $\tilde{Z}$ of
$\mathcal{F}_{\gamma}[\tilde{Z}]=0$, that is, by a suitably scaled
Bryant soliton. These assumptions complete our general \emph{Ansatz}.

\begin{ansatz}
\label{TipAnsatz}
In the tip region, we assume that the \textsc{lhs} of equation~\eqref{Z-time-evolution}
is negligible compared to the \textsc{rhs} for $t \approx T$. To wit, we assume
that (i) $\Gamma(\tau)\gg e^{\tau /2}$, (ii) $\Gamma_\tau=o(\Gamma^{3})$, and
(iii) $Z_{t}=o(\Gamma^{2})$, all as $t\nearrow T$.
\end{ansatz}

As noted above, the choice of the expansion factor $\Gamma(\tau)$ is determined by matching
at the intersection of the tip and parabolic regions. We now discuss this determination. Below,
we verify that the form of $\Gamma(\tau)$ we obtain in equation~\eqref{ChooseGamma}
satisfies the conditions of Ansatz Condition~\ref{TipAnsatz} as long as $k\geq 3$.

To study the consequences of matching at the tip-parabolic intersection, we first
recall that equation \eqref{eigenmode}, with $\sigma=e^{\tau/2}s$, implies that
\begin{equation}\label{match-one}
\psi(s,t)=\sqrt{2(n-1)} e^{(-1/2+\lambda_k)\tau}[1+o(1)]b_k \sigma^{k}
\end{equation}
and
\begin{equation}\label{match-two}
\psi_{s}(s,t)=\sqrt{2(n-1)} e^{\lambda_k \tau}[1+o(1)]k b_k \sigma^{k-1}
\end{equation}
hold for large $|\sigma|$. On the other side of the interface --- at the outer
boundary of the tip region --- it follows from Proposition~\ref{BryantProperties}
in Appendix~\ref{Bryant} and the choice ($c_2=1$) of normalization there that
\begin{equation} \label{match-three}
\psi_{s}(s,t)
    \approx \sqrt{\tilde{Z}(\gamma,t)}
    \approx a\Gamma^{-1}\psi^{-1}(s,t).
\end{equation}

Comparing \eqref{match-one}, \eqref{match-two}, and \eqref{match-three} at $\rho=c$,
i.e.~at $\sigma=e^{-\frac{\lambda_k}{k}\tau}c$, and using \eqref{choose-b}, one finds
that the asymptotic expansions match provided that
\begin{equation}\label{choose-a}
a = \frac{k(n-1)}{2c}
\end{equation}
and
\begin{equation}\label{ChooseGamma}
\Gamma=e^{\left(1-\frac{1}{k}\right)\tau}.
\end{equation}
This yields a one-parameter family of formal solutions indexed by the scaling parameter
$c>0$. More precisely, the free parameters $a$, $b_k$, and $c$ obtained in the tip,
parabolic, and intermediate regions, respectively, are reduced by the tip-parabolic and
parabolic-intermediate matching conditions to a single parameter, which we take to be $c$.

The choice of $\Gamma$ in equation~\eqref{ChooseGamma} satisfies properties~(i) and (ii)
of Ansatz Condition~\ref{TipAnsatz} provided that $k>2$. Moreover, the first-order term
in our approximate solution is stationary, so that property~(iii) is satisfied automatically.
\medskip

We now recall from Section~\ref{Notation} that an $\rm{SO}(n+1)$-invariant
metric has two distinguished sectional curvatures, which we call $K$ and $L$.
In terms of the present set of local coordinates, the sectional curvatures of
$g=z^{-1}(\mathrm{d}\psi)^{2}+\psi^{2}\,g_{\rm{can}}$ are given by
$K=-z_{\psi}/(2\psi)$ and $L=(1-z)/\psi^{2}$. At the pole $x=1$, these
are equal and easily computed. Again by Proposition~\ref{BryantProperties}
in Appendix~\ref{Bryant}, one has
\begin{equation}  \label{BlowUpRate}
K|_{x=1}=L|_{x=1}=\lim_{\psi\searrow 0}\frac{1-Z(\gamma,t)}{\psi^2}
    =\frac{\Gamma^2}{a^2}
    =\frac{a^{-2}}{(T-t)^{2-2/k}},
\end{equation}
where $k\geq 3$. It follows that these formal solutions exhibit the characteristic
Type-II curvature blow-up behavior expected of degenerate neckpinch solutions.
Notice that the limiting behavior as $k\rightarrow\infty$ matches the $(T-t)^{-2}$
blow-up rate first observed by Daskalopoulos and Hamilton in their rigorous treatment
\cite{DH04} of complete Type-II Ricci flow singularities on $\mathbb{R}^2$, where
Ricci flow coincides with the logarithmic fast diffusion equation $u_t=\Delta\log u$.
The asymptotic profile of its blow-up was derived formally by King \cite{King93}
and recovered rigorously by Daskalopoulos and \v{S}e\v{s}um \cite{DS10}.

\section{Conclusions}
\label{TheBitterFinalEnd}

Gu and Zhu \cite{GZ08} prove that Type-II Ricci flow singularities develop from nongeneric
rotationally symmetric Riemannian metrics on $\mathcal{S}^{n+1}\,\,(n\geq 2)$, having the form
$g=(\mathrm{d}s)^{2}+\psi^{2}(s)\,g_{\rm{can}}$ in local coordinates on $\mathcal{S}^{n+1}\backslash\{P_\pm\}$.

Our work above describes and provides plausibility arguments for a detailed asymptotic
profile and rate of curvature blow-up that we predict some (though not necessarily all)
such solutions should exhibit. We summarize our prediction as follows.

\begin{figure}[t]
  \centering
  \includegraphics[width=\textwidth]{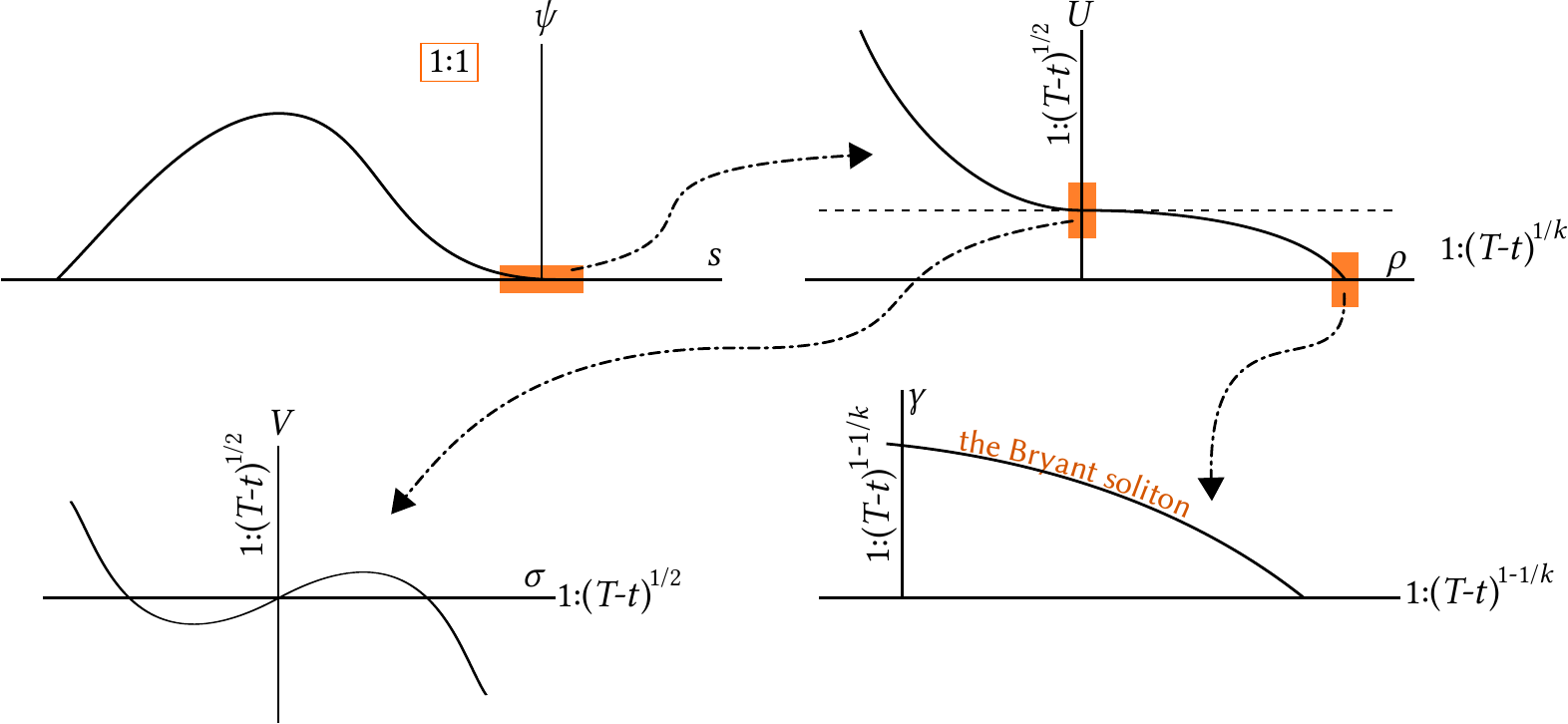}
  \caption{An asymmetric degenerate neckpinch viewed at various scales.}
\end{figure}

\begin{conjecture}

For every $n \geq 2$, every $k\geq 3$, and every $c>0$, there exist Ricci flow
solutions $g(t)$ that satisfy the conditions outlined in our Basic~Assumption and
develop a degenerate neckpinch singularity at the right pole at some $T<\infty$.
The singularity is Type-II --- slowly forming --- with
\[
\sup_{x\in\mathcal{S}^{n+1}}|\mathrm{Rm}(x,t)|\sim\frac{C}{(T-t)^{2-2/k}}
\]
attained at the pole. Its asymptotic profile is as follows, where $s(x,t)$
represents arc-length with respect to $g(t)$ measured from the location
of the smallest nondegenerate neck.

\smallskip
\textsc{Outer Region:} As $t\nearrow T$, one has
\[
\psi(s,t)=[1+o(1)]\sqrt{2(n-1)\left[(T-t)-\left(\frac{s}{c}\right)^k\right]}
\]
holding for $-\varepsilon\leq s\leq c(T-t)^{1/k}$ if $k$ is odd,
and for $|s| \leq c(T-t)^{1/k}$ if $k$ is even.

\smallskip
\textsc{Intermediate Region:} As $t\nearrow T$, one has
\[
\frac{\psi(s,t)}{\sqrt{2(n-1)(T-t)}}=
    [1+o(1)]\sqrt{1-\frac{(s/c)^k}{T-t}}
\]
on an interval
$\varepsilon (T-t)^{1/k}\leq s \leq\varepsilon^{-1}(T-t)^{1/k}$.

\smallskip
\textsc{Parabolic Region:} As $t\nearrow T$, one has
\[
\frac{\psi(s,t)}{\sqrt{2(n-1)(T-t)}}=1-
    \frac{1+o(1)}{2c^k}\frac{\sqrt{(T-t)^{k}}}{T-t}\,\,h_k\left(\frac{s}{\sqrt{T-t}}\right)
\]
on an interval $\varepsilon\sqrt{T-t}\leq s \leq\varepsilon (T-t)^{1/k}$, where
$h_k(\cdot)$ denotes the $k^\mathrm{th}$ Hermite polynomial, normalized so that its
highest-order term has coefficient $1$.

\textsc{Tip Region:} A Bryant soliton forms in a neighborhood of the pole.
Specifically, with respect to a rescaled local radial coordinate
\[
\gamma(s,t)=\frac{\psi(s,t)}{(T-t)^{1-1/k}}
\]
near the right pole, in which the metric takes the form
\[
g=Z(\gamma,t)^{-1}(\mathrm{d}\psi)^2 +\psi^2 \,g_{\rm{can}},
\]
one has
\[
Z(\gamma,t)=[1+o(1)]\,B\left(\frac{2c\gamma}{k(n-1)}\right)\quad (t\nearrow T),
\]
where $B$ denotes the Bryant soliton --- up to scaling, the unique complete
locally conformally flat\,\footnote{Uniqueness also holds if local conformal flatness
is replaced by a suitable condition at spatial infinity; see \cite{Brendle11}.}
non-flat steady gradient soliton on $\mathbb{R}^{n+1}$.
\end{conjecture}
\medskip

\begin{figure}[h]
\centering
\includegraphics[width=\textwidth]{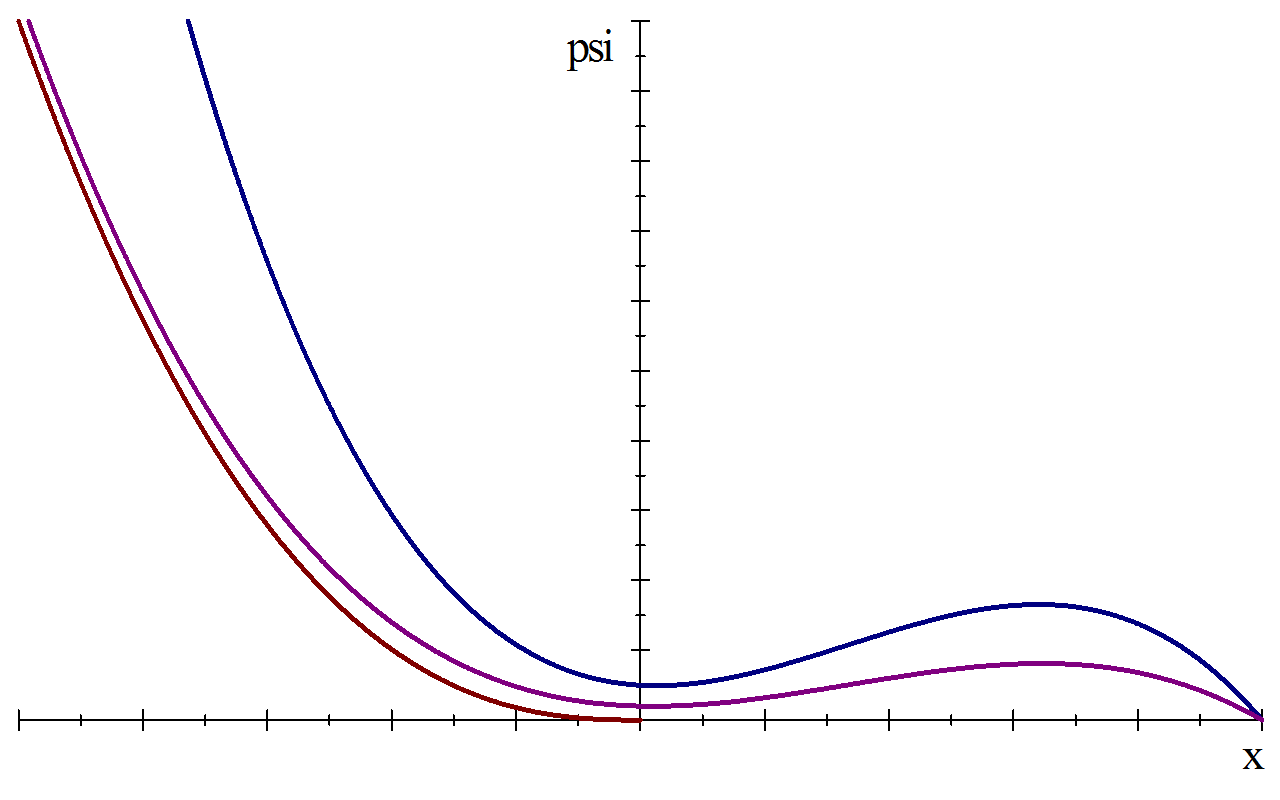}
\caption{An asymmetric degenerate neckpinch viewed without rescaling
in the (non-geometric) $x$ coordinate system.}
\end{figure}

Many aspects of this work are familiar. Indeed, our predicted rate of curvature
blowup matches that of the examples of Type-II mean curvature flow singularities
rigorously constructed by Vel\'{a}zquez and one of the authors \cite{AV97}.
The ``global singularities'' encountered by the symmetric ($k$ even) profiles
considered here agree with the intuition obtained from such rigorous examples
for mean curvature flow. Moreover, the Ricci flow singularities numerically
simulated by Garfinkle and another of the authors \cite{GI03, GI07} are
qualitatively similar to the case $k=4$ considered here.

On the other hand, it was perhaps not obvious \emph{a priori} that a Type-II Ricci flow
solution would vanish on an open set $(0,1)\times\mathcal{S}^{n}$ of the original manifold.
This occurs for the asymmetric ($k$ odd) profiles considered here, which correspond to the
``intuitive solutions'' predicted and sketched by Hamilton \cite[Section~3]{Ham95}.

Motivated by the fact that conjectures provide direction and structure for the development
of rigorous new mathematics, it is our hope that the formal derivations in this paper facilitate
further study of Type-II (degenerate) Ricci flow singularity formation. In particular, we intend
in forthcoming work to provide a rigorous proof that there exist solutions exhibiting the asymptotic
behavior formally described here. Further study is also needed to determine what (if any) other
asymptotic behaviors are possible.

\appendix

\section{Evolution equations in rescaled coordinates}
\label{WhyTheHeck?}
In this appendix, we (partially) explain our choices of space and time dilation in the parabolic
region, and we derive the evolution equations satisfied by a rescaled solution.

Let $T<\infty$ denote the singularity time; let
\begin{equation}
\tau=-\log(T-t);
\end{equation}
and let
\begin{equation}
\sigma=e^{\beta\tau}s,
\end{equation}
where $\beta$ is a constant to be chosen. Define a rescaled solution $U(\sigma,\tau)$ by
\begin{equation}
\psi(s,t)=\sqrt{2(n-1)}e^{-\alpha\tau}U(\sigma,\tau),
\end{equation}
where $\alpha$ is another constant to be determined.

It is straightforward to calculate that
\begin{align*}
\psi_{t}  & =\sqrt{2(n-1)}e^{(1-\alpha)\tau}(U_{\tau}+\sigma_\tau U_{\sigma}-\alpha U),\\
\psi_{s}  & =\sqrt{2(n-1)}e^{(\beta-\alpha)\tau}U_{\sigma},\\
\psi_{ss}  & =\sqrt{2(n-1)}e^{(2\beta-\alpha)\tau}U_{\sigma\sigma}.
\end{align*}
The factor $\sigma_\tau$  is necessary for $\sigma$ and $\tau$ to be commuting variables.
It is given by
\begin{equation}
\sigma_\tau     %%\frac{\partial\sigma}{\partial\tau}
    = \beta\sigma + e^{-\tau/2}\frac{\partial s}{\partial t}
    =  \beta\sigma +nI,
\end{equation}
where $I$ is the non-local term
\begin{equation}\label{non-local}
I := \int_{0}^\sigma\frac{U_{\sigma\sigma}}{U}\,\rm{d}\sigma.
\end{equation}

It follows from equation~\eqref{psi-s-evolution}
that $U$ satisfies
\begin{equation}
e^{(1-2\beta)\tau}(U_{\tau}+\sigma_\tau U_{\sigma}-\alpha U)=U_{\sigma\sigma
}+(n-1)\frac{U_{\sigma}^{2}}{U}-e^{2(\alpha-\beta)\tau}\frac{1}{2U}.
\end{equation}
If $\beta<1/2$, then for $\tau\gg 0$, $U$ should approximate a translating solution of a
first-order equation. If $\beta>1/2$, then for $\tau\gg 0$, $U$ should approximate a
stationary solution of an elliptic equation. The choice $\beta=1/2$ is thus necessary if
one expects $U$ to be modeled by the solution of a parabolic equation.\footnote
    {Caveat: it is not necessarily the case that a parabolic equation will dominate
    in the ``parabolic region;'' for example, see the degenerate singularity considered
    in \cite{ACK09}.}

Now suppose that $\beta=1/2$ and write $U=1+V$. It follows from the considerations above
that $V$ evolves by
\begin{equation}
\label{v-evolution}
V_{\tau}=V_{\sigma\sigma}-\left(\frac{\sigma}{2}+nI\right)V_{\sigma}
    +(n-1)\frac{V_{\sigma}^{2}}{1+V}+\alpha(1+V)-\frac{e^{(2\alpha-1)\tau}}{2(1+V)}.
\end{equation}
Linearizing about $V=0$, i.e.~$U=1$, one finds that
$V_{\tau}=\tilde{A}V+\mathcal{O}(V^{2})$,
where
\[
\tilde{A}:V\mapsto V_{\sigma\sigma}-\frac{\sigma}{2}V_{\sigma}
    +\left[  \alpha+\frac{1}{2}e^{(2\alpha-1)\tau}\right]V
    +\left[  \alpha-\frac{1}{2}e^{(2\alpha-1)\tau}\right].
\]
The choice $\alpha=1/2$ thus results in $\tilde{A}$ becoming the autonomous
linear operator
\begin{equation}
A:V\mapsto V_{\sigma\sigma}-\frac{\sigma}{2}V_{\sigma}+V.
\end{equation}

\section{The Bryant Soliton}
\label{Bryant}

The Bryant soliton, discovered in unpublished work of Robert Bryant, is up to homethetic scaling,
the unique complete non-flat locally conformally flat steady gradient soliton on $\mathbb{R}^{n+1}$
for $n \geq 2$. (See \cite{CC10} and \cite{CM10}.) Uniqueness also holds under the assumption that a 
vector field $V:=\nabla R+\varrho(R)\nabla f$ decays fast enough at spatial infinity. (Recall that $f$ is
the soliton potential function, and $\varrho(R)$ is chosen so that $V$ vanishes on the Bryant soliton;
see \cite{Brendle11}.)

Numerical simulations by Garfinkle and one of the authors suggest that a degenerate neckpinch
solution should converge to the Bryant soliton after rescaling near the north pole \cite{GI03, GI07}.
Results of Gu and Zhu also support this expectation \cite{GZ08}.

Here we recall some relevant properties of these solutions. It is convenient to consider
a one-parameter family of Bryant soliton profile functions $B(\cdot)$ depending on a
positive parameter that encodes the scaling invariance mentioned above. For more
information about the Bryant soliton (including proofs of the following claims), see Appendix~C
of \cite{ACK09} and Chapter~1, Section~4 of \cite{RFV2P1}.

\begin{proposition}[Properties of the Bryant soliton profile function]
\label{BryantProperties}
\quad\\
\begin{enumerate}
\item The \textsc{ode} $\mathcal{F}_r [z]=0$, where
\[
\mathcal{F}_r [z]
    :=\frac{1}{r^2}\left\{r^2 z z_{rr}
    -\frac{1}{2}(r z_r)^2 + (n-1-z)r z_r
    +2(n-1)(1-z)z\right\},
\]
admits a unique one-parameter family of complete solutions satisfying $Z(0)=1$
and $Z(\infty)=0$. These are given by
\[
Z(r)=B\left(\frac{r}{\varrho}\right)
\]
for $\varrho>0$, where $B$
is the \emph{Bryant soliton profile function}. Each member of the one-parameter
family of complete smooth metrics given by
\[
g=Z^{-1}(r)(\mathrm{d}r)^2 +r^2 \,g_{\rm{can}},
\]
is called a \emph{Bryant soliton}.

\item $B(r)$ is strictly monotone decreasing for all $r>0$.

\item Near $r=0$, $B$ is smooth and has the asymptotic expansion
\begin{equation*}
B(r)=1+b_{2}r^{2}+\frac{n}{n+3}b_{2}^{2}r^{4}+\frac{n(n-1)}{(n+3)(n+5)}b_{2}^{3}r^{6}+\cdots,
\end{equation*}
where $b_{2}<0$ is arbitrary.

\item Near $r=+\infty $, $B$ is smooth and has the asymptotic expansion
\begin{equation*}
B(r)=c_{2}r^{-2}+\frac{4-n}{n-1}c_{2}^{2}r^{-4}+\frac{(n-4)(n-7)}{(n-1)^{2}}c_{2}^{3}r^{-6}+\cdots,
\end{equation*}%
where $c_{2}>0$ is arbitrary.
\end{enumerate}
\end{proposition}

The arbitrariness of $b_{2}$ and $c_{2}$ encodes the scaling invariance of
the Bryant soliton. In this paper, we fix $c_2=1$ in order to make
explicit the dependence on the scaling parameter $\varrho>0$ when writing
 $Z(r)=B(r/\varrho)$.


\begin{thebibliography}{9}
\bibitem{ACK09} \textbf{Angenent, Sigurd B.; Caputo, M.~Cristina; Knopf, Dan.}
Minimally invasive surgery for Ricci flow singularities.
\emph{J.~Reine Angew.~Math.} To appear. (\texttt{arXiv:0907.0232})

\bibitem{AK04} \textbf{Angenent, Sigurd B.; Knopf, Dan.}
An example of neckpinching for Ricci flow on $S^{n+1}$.
\emph{Math.~Res.~Lett. } \textbf{11 } (2004), no.~4, 493--518.

\bibitem{AK07} \textbf{Angenent, Sigurd B.; Knopf, Dan.}
Precise asymptotics of the Ricci flow neckpinch.
\emph{Comm.~Anal.~Geom.~}\textbf{15} (2007), no.~4, 773--844.

\bibitem{AV97} \textbf{Angenent, Sigurd B.; Vel\'{a}zquez, J. J. L.}
Degenerate neckpinches in mean curvature flow.
\emph{J.~Reine Angew.~Math.}~\textbf{482} (1997), 15--66.

\bibitem{Brendle11} \textbf{Brendle, Simon.}
Uniqueness of gradient Ricci solitons. \texttt{arXiv:1010.3684}

\bibitem{CC10} \textbf{Cao, Huai-Dong; Chen, Qiang.}
On Locally Conformally Flat Gradient Steady Ricci Solitons. \texttt{arXiv:0909.2833}

\bibitem{CM10} \textbf{Catino, Giovanni; Mantegazza, Carlo.}
Evolution of the Weyl tensor under the Ricci flow.
\emph{Ann.~Inst.~Fourier.} To appear. (\texttt{arXiv:0910.4761v6})

\bibitem{RFV2P1} \textbf{Chow, Bennett; Chu, Sun-Chin; Glickenstein, David;
Guenther, Christine; Isenberg, James; Knopf, Dan; Ivey, Tom; Lu, Peng; Luo, Feng; Ni, Lei.}
\emph{The Ricci Flow: Techniques and Applications, Part I: Geometric Aspects.}
Mathematical Surveys and Monographs, Vol.~135.~American Mathematical Society, Providence, RI, 2007.

\bibitem{DH04} \textbf{Daskalopoulos, Panagiota; Hamilton, Richard S.}
Geometric estimates for the logarithmic fast diffusion equation.
\emph{Comm.~Anal.~Geom.~}\textbf{12} (2004), no.~1-2, 143--164.

\bibitem{DS10} \textbf{Daskalopoulos, Panagiota; \v{S}e\v{s}um, Nata\v{s}a.}
Type II extinction profile of maximal solutions to the Ricci flow in $\Bbb R^2$.
\emph{J.~Geom.~Anal.}~\textbf{20} (2010), no.~3, 565--591.

\bibitem{EMT10} \textbf{Enders, Joerg; M\"{u}ller, Reto; Topping, Peter M.}
On Type I Singularities in Ricci flow. \texttt{arXiv:1005.1624}

\bibitem{FIK03} \textbf{Feldman, Mikhail; Ilmanen, Tom; Knopf, Dan.}
Rotationally symmetric shrinking and expanding gradient K\"{a}hler--Ricci solitons.
\emph{J.~Differential Geom.~}\textbf{65} (2003), no.~2, 169--209.

\bibitem{GI03} \textbf{Garfinkle, David; Isenberg, James.}
Critical behavior in Ricci flow. \texttt{arXiv:math/0306129}.

\bibitem{GI07} \textbf{Garfinkle, David; Isenberg, James.}
The Modelling of Degenerate Neck Pinch Singularities in Ricci Flow by Bryant Solitons.
\texttt{arXiv:0709.0514}.

\bibitem{GK85} \textbf{Giga, Yoshikazu; Kohn, Robert V.}
Asymptotically self-similar blow-up of semilinear heat equations.
\emph{Comm.~Pure Appl.~Math.}~\textbf{38} (1985), no.~3, 297--319.

\bibitem{GZ08} \textbf{Gu, Hui-Ling; Zhu, Xi-Ping.}
The Existence of Type II Singularities for the Ricci Flow on $S^{n+1}$.
\emph{Comm.~Anal.~Geom.}~\textbf{16} (2008), no.~3, 467--494.

\bibitem {Ham95} \textbf{Hamilton, Richard S.}
The formation of singularities in the Ricci flow.\ \emph{Surveys in
differential geometry, Vol.~II} (Cambridge, MA, 1993), 7--136,
Internat.~Press, Cambridge, MA, 1995.

\bibitem{King93} \textbf{King, John R.}
Self-similar behavior for the equation of fast nonlinear diffusion.
\emph{Philos.~Trans.~R.~Soc., Lond.}, A, \textbf{343}, (1993) 337--375.

\bibitem {KL08} \textbf{Kleiner, Bruce; Lott, John.}
Notes on Perelman's papers.
\emph{Geom.Topol.}~\textbf{12} (2008), no.~5, 2587--2855.

\bibitem{P1} \textbf{Perelman, Grisha.}
The entropy formula for the Ricci flow and its geometric applications.
\texttt{arXiv:math.DG/0211159}.

\bibitem{Simon00}\textbf{Simon, Miles}.
A class of Riemannian manifolds that pinch when evolved by Ricci flow.
\emph{Manuscripta Math.~}\textbf{101} (2000), no.~1, 89--114.

\end{thebibliography}
\end{document}